\documentclass[11pt]{article}
\usepackage{amsmath}
\usepackage{amssymb}
\usepackage{amsthm}
\usepackage[usenames]{color}
\usepackage{colortbl}
\def\R{\mathbb{R}}
\def\N{\mathbb{N}}
\def\C{\mathbb{C}}
\def\a{\mathbf{a}}
\def\b{\mathbf{b}}
\def\c{\mathbf{c}}
\def\d{\mathbf{d}}

\newtheorem{theorem}{\hspace*{\parindent}Theorem}
\newtheorem{lemma}{\hspace*{\parindent}Lemma}
\newtheorem{corollary}{\hspace*{\parindent}Corollary}

\newcounter{theremark}

\title{Inequalities for  some basic hypergeometric functions}
\title{Inequalities for  some basic hypergeometric functions}
\author{S.I.\:Kalmykov$^{\rm a,b}$ and D.B.\:Karp$^{\rm a,b}$\footnote{Corresponding author. E-mail: S.I.\:Kalmykov -- \emph{sergeykalmykov@inbox.ru}, D.B.\:Karp -- \emph{dimkrp@gmail.com}}
\\[10pt]
\small{\textit{$\phantom{1}^a$Far Eastern Federal University, Vladivostok, Russia}}
\\\small{\textit{$\phantom{1}^b$Institute of Applied Mathematics, FEBRAS, Vladivostok, Russia}}}
\date{}

\begin{document}

\maketitle

\begin{abstract}
We establish conditions for the discrete versions of logarithmic concavity and convexity of the higher order regularized basic hypergeometric function with respect
simultaneous shift of all its parameters.  For a particular case of Heine's basic hypergeometric function we prove logarithmic concavity and convexity with respect to the bottom parameter.  We further establish a linearization identity for the generalized Tur\'{a}nian  formed by a particular case of Heine's basic hypergeometric function.  Its $q=1$ case also appears to be new.
\end{abstract}

\bigskip


\bigskip


\section{Introduction}

We will use the standard definition of the $q$-shifted factorial \cite[(1.2.15)]{GR}:
\begin{equation}\label{eq:q-shifted}
(a;q)_0=1,~~~~(a;q)_n=\prod\limits_{k=0}^{n-1}(1-aq^k),~~n\in\N.
\end{equation}
This definition works for any complex $a$ and $q$ but in this paper we confine ourselves to the case $0<q<1$.
Under this restriction we can also define
$$
(a;q)_{\infty}=\lim\limits_{n\to\infty}(a;q)_n,
$$
where the limit can be shown to exist as a finite number for all complex $a$.  The $q$-gamma function is given by  \cite[(1.10.1)]{GR}, \cite[(21.16)]{KacCheung}
\begin{equation}\label{eq:q-Gamma-def}
\Gamma_q(z)=(1-q)^{1-z}\frac{(q;q)_{\infty}}{(q^z;q)_{\infty}}
\end{equation}
for $|q|<1$ and all complex $z$ such that $q^{z+k}\ne1$ for $k\in\N_0$.

In our recent paper \cite{KKJMAA2018} we studied logarithmic concavity and convexity for generic series
with respect to a parameter contained in the argument of $q$-shifted factorial or $q$-gamma function.  Namely, we have considered the series of the form
\begin{equation}\label{eq:f-general}
f(\mu;x)=\sum\nolimits_{k=0}^{\infty}f_k\phi_k(\mu)x^k,
\end{equation}
where $\phi_k(\mu)$ is one of the functions $(q^{\mu};q)_{k}$, $[(q^{\mu};q)_k]^{-1}$, $\Gamma_q(\mu+k)$ or $\Gamma_q^{-1}(\mu+k)$, and $f_k$ is a (typically nonnegative) real sequence independent of $\mu$ and $x$.  Our results were given in terms of the sign of the generalized Tur\'{a}nian
\begin{equation}\label{eq:delta-defined}
\Delta_f(\alpha,\beta;x)=f(\mu+\alpha;x)f(\mu+\beta;x)-f(\mu;x)f(\mu+\alpha+\beta;x)=\sum\limits_{m=0}^{\infty}\delta_mx^m
\end{equation}
 and  its power series coefficients $\delta_m$. Note that $\Delta_f(\alpha,\beta;x)\ge0$ for
$\mu$, $\mu+\alpha$, $\mu+\beta$ and $\mu+\alpha+\beta$ in some interval and $\alpha,\beta\ge0$ is equivalent to the logarithmic concavity of $\mu\to f(\mu;x)$ in that interval. Natural examples of the series (\ref{eq:f-general}) come from the realm of the basic (or $q$-) hypergeometric functions: seven such examples were discussed by us in \cite{KKJMAA2018}.

The purpose of this paper is to prove  new inequalities for the basic hypergeometric functions most of which cannot be derived from the general theorems presented in \cite{KKJMAA2018}.  Furthermore, we establish several related identities some of which are believed to be new. We remark here that while identities for the basic hypergeometric functions form
an immense subject leaving no hope for making any reasonable survey of, inequalities for and convexity-like properties of these functions seem to be a rather rare species.  We are only aware of a handful of papers dealing with this topic. Other than our paper \cite{KKJMAA2018}, these include the works of Baricz, Raghavendar and Swaminathan \cite{BRS}, Mehrez \cite{Mehrez} and Mehrez and Sitnik \cite{MehrSitn}.   An inequality for the basic hypergeometric function of a different type was also established by Zhang in \cite{Zhang}.

\section{Results for the generalized $q$-hypergeometric series}

Throughout the paper we will use the following short-hand notation for products and sums. Given a vector $\a=(a_1,\ldots,a_m)$ and a scalar $\mu$, write
$$
\a+\mu=(a_1+\mu,\ldots,a_m+\mu),~~q^{\a}=q^{a_1}\cdots{q^{a_m}},~~(1-q^{\a})=(1-q^{a_1})\cdots(1-q^{a_m}),
$$
$$
(\a; q)_{n}=(a_1,\ldots,a_m;q)_n=(a_1;q)_{n}(a_2;q)_{n}\cdots(a_m;q)_{n},~~~\Gamma_q(\a)=\Gamma_q(a_1)\cdots\Gamma_q(a_m),
$$
where $n$ may take the value $\infty$.  The generalized $q$-hypergeometric series is defined by \cite[formula (1.2.22)]{GR}
\begin{equation}\label{eq:r-phi-s}
{_t\phi_s}\left(\!\begin{array}{l}\a\\\b\end{array}\!\vline\,q;z\right)={_t\phi_s}(\a;\b;q;z)=\sum_{n=0}^{\infty}\frac{(\a;q)_n}{(\b;q)_n(q;q)_n} \left[(-1)^nq^{\binom{n}{2}}\right]^{1+s-t}z^n,
\end{equation}
where $\a=(a_1,\ldots, a_t)$, $\b=(b_1,\ldots,b_s)$, $t\le{s+1}$, and the series converges for all $z$ if $t\le{s}$ and for $|z|<1$ if $t=s+1$ \cite[section~1.2]{GR}.  As the base $q$ will remain fixed throughout the paper we will omit $q$ from the above notation and write ${_t\phi_s}(\a;\b;z)$ for the right hand side of (\ref{eq:r-phi-s}).

In this section we will deal with the function
\begin{equation}\label{eq:genqnorm}
g(\mu;x)=\frac{\Gamma_q(\a+\mu)}{\Gamma_q(\b+\mu)}
{}_t\phi_s\left(\!\begin{array}{c} q^{\a+\mu}\\q^{\b+\mu}\!\end{array}\vline\,(q-1)^{1+s-t}x\right),
\end{equation}
where ${}_t\phi_s$ is defined in (\ref{eq:r-phi-s}).  The main result will be formulated in terms of the generalized Tur\'anian
\begin{equation}\label{eq:genturqnorm}
\Delta_g(\alpha,\beta;x):=g(\mu+\alpha;x)g(\mu+\beta;x)-g(\mu;x)g(\mu+\alpha+\beta;x)=\sum_{m=0}^{\infty}\gamma_mx^m.
\end{equation}
We will give sufficient conditions for its positivity and negativity as well as for its power series coefficients $\gamma_m$. To this end, we need the standard definition of the elementary symmetric polynomials \cite[3.F]{MOA}:
$$
e_k(c_1,\dots,c_r)=\sum\limits_{1\le{j_1}<j_2\cdots<j_k\le{r}}c_{j_1}c_{j_2}\cdots{c_{j_k}},~~~k=1,\ldots,r,
$$
and $e_0(c_1,\dots,c_r)=1$. The following theorem is a $q$-analogue of \cite[Theorem~3]{KKJA2016}.
\begin{theorem}\label{th:tphis}
For a given $0<q<1$ and nonnegative vectors $\a$ and $\b$,  set $\c=(q^{-a_1}-1,\dots, q^{-a_t}-1)$ and $\d=(q^{-b_1}-1,\dots,q^{-b_{s}}-1)$.  The following statements hold for $\Delta_g(\alpha,\beta;x)$ defined in \emph{(\ref{eq:genturqnorm}):}

\medskip

\emph{(a)} If $s\le t\le s+1$ and the conditions
\begin{equation}\label{eq:ratincr}
\frac{e_t(\c)}{e_s(\d)}\leq \frac{e_{t-1}(\c)}{e_{s-1}(\d)}\leq \dots \leq \frac{e_{t-s+1}(\c)}{e_1(\d)}\leq e_{t-s}(\c)
\end{equation}
are satisfied, then   $\Delta_g(\alpha,\beta;x)\le 0$ for  $x\ge0$, $\beta\ge0$ and $\alpha\in\N$. Moreover, if $\alpha\le\beta+1$, then the power series coefficients $\gamma_m$ of $x\to\Delta_g(\alpha,\beta;x)$ are non-positive. In particular, the inverse Tur\'{a}n type inequality $g(\mu+1;x)^2\le g(\mu;x)g(\mu+2;x)$ holds for $\mu\ge0$ and fixed $x>0$ in the domain of convergence.

\medskip

\emph{(b)} If $t\le s$ and the conditions
\begin{equation}\label{eq:ratdecr}
\frac{e_s(\d)}{e_t(\c)}\leq \frac{e_{s-1}(\d)}{e_{t-1}(\c)}\leq\cdots\leq \frac{e_{s-t+1}(\d)}{e_1(\c)}\leq e_{s-t}(\d)
\end{equation}
are satisfied, then   $\Delta_g(\alpha,\beta;x)\ge 0$ for  $x\ge0$, $\beta\ge0$ and $\alpha\in\N$. Moreover, if $\alpha\le\beta+1$, then the power series coefficients $\gamma_m$ of $x\to\Delta_f(\alpha,\beta;x)$ are non-negative. In particular, the Tur\'{a}n type inequality
$g(\mu+1;x)^2\ge g(\mu;x)g(\mu+2;x)$ holds for $\mu\ge0$ and any fixed $x>0$.
\end{theorem}

In many situations complicated conditions (\ref{eq:ratincr}) or (\ref{eq:ratdecr}) can be replaced by slightly stronger but simpler majorization conditions: for given vectors $\c,\d\in\R^t$ the weak supermajorization $\d\prec^W\c$ means that \cite[Definition~A.2]{MOA}
\begin{equation*}\label{eq:amajorb}
\begin{split}
& 0<c_1\leq{c_2}\leq\cdots\leq{c_t},~~
0<d_1\leq{d_2}\leq\cdots\leq{d_t},
\\
&\sum\limits_{i=1}^{k}c_i\leq\sum\limits_{i=1}^{k}d_i~~\text{for}~~k=1,2\ldots,t.
\end{split}
\end{equation*}

Denote the size of a vector $\a$ by $|\a|$. The following result was demonstrated by us in \cite[Lemma~4]{KKJA2016}.
\begin{lemma}\label{lm:majincdecr}
Let $\c\in\R^t$, $\d\in\R^s$ be positive vectors, $t\ge{s}$, and suppose that there exists $\c'\subset\c$, $|\c'|=s$, such that $\d\prec^W\c'$.  Then inequalities \emph{(\ref{eq:ratincr})} hold. Similarly, inequalities \emph{(\ref{eq:ratdecr})} hold if $t\le{s}$ and $\c\prec^W\d'$ for some subvector $\d'\subset\d$ of size $t$.
\end{lemma}

Recall that a nonnegative function $f$ defined on an interval $I$ is called completely monotonic there if it has derivatives of all orders and $(-1)^nf^{(n)}(x)\ge0$ for $n\in\N_0$ and $x\in{I}$, see \cite[Defintion~1.3]{SSV}.   We need the following simple lemma.  See \cite[Lemma~2]{KKJMAA2018} and \cite[Remark, p.340]{KKJMAA2018} for a proof.
\begin{lemma}\label{lm:multconvex}
Suppose $\phi(x)=\sum_{k\ge0}\phi_kx^k$ converges for $|x|<R$ with $0<R\le\infty$ and $\phi_k\ge0$. Then $x\to\phi(x)$ is multiplicatively convex and $y\to\phi(1/y)$ is completely monotonic \emph{(}and hence is log-convex\emph{)} on $(1/R,\infty)$, so that there exists a non-negative measure $\tau$ supported on $[0,\infty)$ such that
$$
\phi(x)=\int_{[0,\infty)}e^{-(1/x-1/R)t}\tau(dt).
$$
If $R=\infty$, this measure is given by
$$
\tau(dt)=\phi_0\mathbf{1}_{0}+\left(\sum\nolimits_{m=1}^{\infty}\frac{\phi_mt^m}{(m-1)!}\right)dt,
$$
where $\mathbf{1}_{0}$ is the unit mass concentrated at zero.
\end{lemma}

We further remark that for $R=\infty$ the function $\phi(1/y)$ satisfies the conditions of \cite[Theorem~1.1]{KoumPeders}
and hence enjoys all the properties stated in that theorem. Application of Lemma~\ref{lm:multconvex} leads immediately to the following corollary to Theorem~\ref{th:tphis}.

\begin{corollary}\label{cor:repr}
Suppose conditions of Theorem~\emph{\ref{th:tphis}(b)} are satisfied. Then the $x\mapsto\Delta_g(\alpha,\beta;x)$ defined in \emph{(\ref{eq:genturqnorm})} is multiplicatively convex on $(0,\infty)$, while the function $y\mapsto\Delta_g(\alpha,\beta;1/y)$ is completely monotonic \emph{(}and therefore log-convex\emph{)} on $(0,\infty)$, so that
$$
\Delta_g(\alpha,\beta;x)=\int_{[0,\infty)}e^{-t/x}\tau(dt),
$$
where the non-negative measure $\tau$ is given by
$$
\tau(dt)=\gamma_0\mathbf{1}_{0}+\left(\sum\nolimits_{m=1}^{\infty}\frac{\gamma_mt^m}{(m-1)!}\right)dt.
$$
Here $\mathbf{1}_{0}$ stands for the unit mass concentrated at zero and $\gamma_m$ are defined in \emph{(\ref{eq:genturqnorm})}.

If conditions of Theorem~\emph{\ref{th:tphis}(a)} are satisfied and $t=s$ similar statements hold for $-\Delta_g(\alpha,\beta;x)$.
\end{corollary}

We conclude this section with two examples of the application of Theorem~\ref{th:tphis}.

\medskip

\noindent\textbf{Example~1.}
Suppose $0<a_1\le{a_2}$, $0<b_1\le{b_2}$ and
\begin{equation}\label{eq:alphabetacond}
\begin{split}
& q^{-b_1}\leq q^{-a_1}~\Leftrightarrow~b_1\le{a_1},
\\
& q^{-b_1}+q^{-b_2}\leq q^{-a_1}+q^{-a_2}.
\end{split}
\end{equation}
Then by Lemma~\ref{lm:majincdecr} and Theorem~\ref{th:tphis}(b) (with $t=s=2$) the generalized Tur\'{a}nian $\Delta_g(\alpha,\beta;x)$ with
$$
g(\mu;x)=\frac{\Gamma_q(a_1+\mu)\Gamma_q(a_2+\mu)}{\Gamma_q(b_1+\mu)\Gamma_q(b_2+\mu)}
{}_2\phi_2\left(\!\begin{array}{c} q^{a_1+\mu}, q^{a_2+\mu}\\q^{b_1+\mu}, q^{b_2+\mu}\!\end{array}\vline\,(q-1)x\right)
$$
is non-negative for  $x\ge0$, $\beta\ge0$ and $\alpha\in\N$. Moreover, if $\alpha\le\beta+1$, then the power series coefficients of $x\to\Delta_g(\alpha,\beta;x)$ are all non-negative.  If both inequalities (\ref{eq:alphabetacond}) are reversed, then $x\to\Delta_g(\alpha,\beta;x)$ has non-positive power series coefficients for  $x\ge0$, $\beta\ge0$ and $\N\ni\alpha\le\beta+1$ and $\Delta_g(\alpha,\beta;x)\le0$ holds without the restriction $\alpha\le\beta+1$.

\bigskip

\noindent\textbf{Example~2.}
Suppose $0<a_1\le{a_2}\le{a_3}$, $0<b_1\le{b_2}$ and
\begin{equation*}
\begin{split}
& q^{-a_1}\leq q^{-b_1}~\Leftrightarrow~a_1\le{b_1},
\\
& q^{-a_1}+q^{-a_2}\leq q^{-b_1}+q^{-b_2}.
\end{split}
\end{equation*}
Then by Lemma~\ref{lm:majincdecr} and Theorem~\ref{th:tphis}(a) (with $t=3$, $s=2$) the generalized Tur\'{a}nian $\Delta_g(\alpha,\beta;x)$ with
$$
g(\mu;x)=\frac{\Gamma_q(a_1+\mu)\Gamma_q(a_2+\mu)\Gamma_q(a_3+\mu)}{\Gamma_q(b_1+\mu)\Gamma_q(b_2+\mu)}
{}_3\phi_2\left(\!\begin{array}{c} q^{a_1+\mu}, q^{a_2+\mu}, q^{a_3+\mu}\\q^{b_1+\mu}, q^{b_2+\mu}\!\end{array}\vline\,x\right)
$$
is non-positive for  $x\ge0$, $\beta\ge0$ and $\alpha\in\N$. Moreover, if $\alpha\le\beta+1$, then the power series coefficients of $x\to\Delta_g(\alpha,\beta;x)$ are all non-positive.

\section{Results for Heine's $q$-hypergeometric series}

Before stating the results we need some definitions.  The (first) $q$-exponential is defined by  \cite[formula~(18)]{GR}
\begin{equation}\label{eq:symmetricpol}
e_q(z)=\sum\limits_{k=0}^{\infty}\frac{z^k}{(q;q)_k},~~~~~~|z|<1.
\end{equation}
Jackson gave the following two $q$-analogues of the Bessel functions \cite[Exersice ~1.24, p.30]{GR}:
\begin{equation}\label{eq:J1}
J^{(1)}_{\alpha}(y)=\frac{(y/2)^{\alpha}(q^{\alpha+1};q)_{\infty}}{(q;q)_{\infty}}{_2\phi_1}(0,0;q^{\alpha+1};-y^2/4),~~~~|y|<2,
\end{equation}
\begin{equation}\label{eq:J2}
J^{(2)}_{\alpha}(y)=\frac{(y/2)^{\alpha}(q^{\alpha+1};q)_{\infty}}{(q;q)_{\infty}}{_0\phi_1}(-;q^{\alpha+1};-y^2q^{\alpha+1}/4),~~~y\in\C.
\end{equation}

The lemma below is a slight modification of a product formula due to Rahman \cite[formula (1.20)]{Rahman}.  We present it in this section as it may be of independent interest.
\begin{lemma}\label{lm:qBessel}
The following $q$-identity holds when both sides are well defined
\begin{equation}\label{eq:qBessel}
{_2\phi_1}(0,0;q^{\nu};z){_2\phi_1}(0,0;q^{\eta};z)=e_q(z){_4\phi_3}\!\!
\left(\left.\!\!\!\!\begin{array}{c}q^{(\nu+\eta-1)/2},q^{(\nu+\eta)/2},-q^{(\nu+\eta-1)/2},-q^{(\nu+\eta)/2}\\
q^{\nu},q^{\eta},q^{\nu+\eta-1}\end{array}\right|z\!\right).
\end{equation}
\end{lemma}
\textbf{Proof.} The identity given in \cite[formula (1.20)]{Rahman} can be written in terms of the second $q$-Bessel function $J^{(2)}_{\alpha}$ as follows:
\begin{multline*}
J^{(2)}_{\alpha}(y)J^{(2)}_{\beta}(y)/(-y^2/4;q)_{\infty}=
\frac{y^{\alpha+\beta}}{(2(1-q))^{\alpha+\beta}\Gamma_q(\alpha+1)\Gamma_q(\beta+1)}\times
\\
{_4\phi_3}\!\!
\left(\left.\!\!\!\!\begin{array}{c}q^{(\alpha+\beta+1)/2},q^{(\alpha+\beta+2)/2},-q^{(\alpha+\beta+1)/2},-q^{(\alpha+\beta+2)/2}\\
q^{\alpha+1},q^{\beta+1},q^{\alpha+\beta+1}\end{array}\right|-\frac{y^2}{4}\!\right).
\end{multline*}
Using the definitions (\ref{eq:J1}), (\ref{eq:J2}) and the connection formula  \cite[Exercise~1.24]{GR}
$$
J^{(2)}_{\alpha}(y)=(-y^2/4;q)_{\infty}J^{(1)}_{\alpha}(y)
$$
we arrive at (\ref{eq:qBessel}) on writing $\nu=\alpha+1$, $\eta=\beta+1$, $z=-y^4/4$ and applying
the summation formula  \cite[(1.14)]{Rahman}, \cite[(1.3.15)]{GR}
$$
e_q(z)=\frac{1}{(z;q)_{\infty}}. ~~~~~~~~~~~~~~\square
$$
\medskip

Equating the coefficients at equal powers of $z$ in (\ref{eq:qBessel}) we immediately obtain
\begin{corollary}\label{cr:qBessel}
For arbitrary indeterminates $\nu$, $\eta$ and $m\in\N$ the identity
\begin{equation}\label{eq:qBesselfinitesum}
\sum\limits_{k=0}^{m}\frac{1}{(q^{\nu};q)_k(q^{\eta};q)_{m-k}(q;q)_{k}(q;q)_{m-k}}
=\sum\limits_{k=0}^{m}\frac{(q^{(\nu+\eta-1)/2},q^{(\nu+\eta)/2},-q^{(\nu+\eta-1)/2},-q^{(\nu+\eta)/2};q)_k}
{(q^{\nu},q^{\eta},q^{\nu+\eta-1};q)_k(q;q)_{k}(q;q)_{m-k}}
\end{equation}
holds true.
\end{corollary}

The following theorem can be established by an application of \cite[Theorem~4]{KKJMAA2018}.  Here we will furnish a simple and independent direct proof based on Lemma~\ref{lm:qBessel}.
\begin{theorem}\label{th:1phi1}
Suppose
$$
f(\mu;x)={}_2\phi_{1}(0,0;q^{\mu};x)=\sum\limits_{n=0}^{\infty}\frac{x^n}{(q^{\mu};q)_n(q;q)_n}.
$$
Then the power series coefficients of the generalized Tur\'{a}nian
$$
\Delta_f(\alpha,\beta;x)=f(\mu+\alpha;x)f(\mu+\beta;x)-f(\mu;x)f(\mu+\alpha+\beta;x)
$$
are negative for all $\mu,\alpha,\beta>0$.  In particular, $\mu\to f(\mu;x)$ is log-convex on $[0,\infty)$ for  $0<x<1$.
\end{theorem}

\noindent\textbf{Remark.} The function appearing in the above theorem can also be written in terms of the (first) modified $q$-Bessel function introduced
by Ismail \cite[(2.5)]{Ism} and rediscovered by Olshanetskii and Rogov \cite[section~3.1]{OR1} as follows:
\begin{equation}\label{eq:2phi1Bessel}
{}_2\phi_{1}(0,0;q^{\mu};x)=\frac{(1-q)^{\mu-1}\Gamma_q(\mu)}{x^{(\mu-1)/2}}I^{(1)}_{\mu-1}(2\sqrt{x};q).
\end{equation}

The basic hypergeometric series treated in the theorem below is a particular case of the generic series considered in \cite[Theorem~3]{KKJMAA2018}.  However, the conclusion we draw here is stronger than that of \cite[Theorem~3]{KKJMAA2018}, where an additional restriction $\alpha\in\N$, $\alpha\le\beta+1$ is imposed.

\begin{theorem}\label{th:2phi1logconcavity}
Suppose
$$
\tilde{f}(\mu;x)=\frac{{}_2\phi_{1}(0,0;q^{\mu};x)}{\Gamma_q(\mu)}
=\sum\limits_{n=0}^{\infty}\frac{[x/(1-q)]^n}{\Gamma_q(\mu+n)(q;q)_n}.
$$
Then the power series coefficients of the generalized Tur\'{a}nian
$$
\Delta_{\tilde{f}}(\alpha,\beta;x)=\tilde{f}(\mu+\alpha;x)\tilde{f}(\mu+\beta;x)-\tilde{f}(\mu;x)\tilde{f}(\mu+\alpha+\beta;x)
$$
are positive for $\mu,\alpha,\beta>0$. In particular, the function $\mu\to\tilde{f}(\mu;x)$ is log-concave on $[0,\infty)$ for $0<x<1$.
\end{theorem}

In  \cite[Example~1]{KKJMAA2018} we announced without proof that the modified $q$-Bessel function $I^{(1)}_{\nu}(y;q)$ is log-concave with respect to $\nu$.
In view of (\ref{eq:2phi1Bessel}), Theorem~\ref{th:2phi1logconcavity} immediately leads to
\begin{corollary}\label{cr:BesslIq}
For each fixed $y\in(0,2)$ the function
$$
\nu\to I^{(1)}_{\nu}(y;q)
$$
is log-concave on $(-1,\infty)$.
\end{corollary}

Finally, we present a  linearization identity for product difference of Heine's $q$-hypergeometric functions ${}_2\phi_1$ which appears to be new.
\begin{theorem}\label{th:2phi1-identity}
For $\alpha\in\N$ the identity
\begin{multline}\label{eq:2phi1-identity}
(q^{\mu+\beta};q)_{\alpha}{}_2\phi_1\left(\!\begin{array}{c}q, 0\\q^{\mu+\alpha}\end{array}\vline\,x\right)
{}_2\phi_1\left(\!\begin{array}{c}q, 0\\q^{\mu+\beta}\end{array}\vline\,x\right)
-(q^{\mu};q)_{\alpha}{}_2\phi_1\left(\!\begin{array}{c}q,0\\q^{\mu}\!\end{array}\vline\,x\right)
{}_2\phi_1\left(\!\begin{array}{c}q,0\\q^{\mu+\alpha+\beta}\!\end{array}\vline\,\,x\right)
\\
= \sum_{j=0}^{\alpha-1}\biggl\{(q^{\mu+1+j};q)_{\alpha-1-j}
(q^{\mu+\alpha+\beta-1-j};q)_{1+j}{}_2\phi_1\left(\!\begin{array}{c}q,0\\q^{\mu+1+j}\!\end{array}\vline\,x\right)
\\
-(q^{\mu+j};q)_{\alpha-j}(q^{\mu+\alpha+\beta-j};q)_{j}
{}_2\phi_1\left(\!\begin{array}{c}q,0\\q^{\mu+\alpha+\beta-j}\!\end{array}\vline\,x\right)\biggr\}
\end{multline}
holds for values of parameters for which both sides make sense.
\end{theorem}

Taking the limit $q\to1$ in the above theorem we arrive at an identity generalizing \cite[Theorem~3]{KKSEMR2018}.
\begin{corollary}\label{cr:Kummer}
For $\alpha\in \N$ the identity
\begin{multline}\label{eq:1F1-identity}
(\mu+\beta)_{\alpha}{}_1F_1\left(\!\begin{array}{c}1 \\ \mu+\alpha\end{array}\vline\,x\right)
{}_1F_1\left(\!\begin{array}{c}1\\ \mu+\beta \end{array}\vline\,x\right)
-(\mu)_{\alpha}{}_1F_1\left(\!\begin{array}{c}1\\ \mu\!\end{array}\vline\,x\right)
{}_1F_1\left(\!\begin{array}{c}1 \\ \mu+\alpha+\beta\!\end{array}\vline\,\,x\right)
\\
= \sum_{j=0}^{\alpha-1}\biggl\{(\mu+1+j)_{\alpha-1-j}
(\mu+\alpha+\beta-1-j)_{1+j}{}_1F_1\left(\!\begin{array}{c}1\\ \mu+1+j\!\end{array}\vline\,x\right)
\\
-(\mu+j)_{\alpha-j}(\mu+\alpha+\beta-j)_{j}
{}_1F_1\left(\!\begin{array}{c}1\\ \mu+\alpha+\beta-j\!\end{array}\vline\,x\right)\biggr\}
\end{multline}
holds for values of parameters for which both sides make sense. Here $(\mu)_{\alpha}$ is the Pochhammer symbol and ${}_1F_1$ stands for the Kummer (or confluent hypergeometric) function \emph{\cite[(1.2.16)]{GR}}.
\end{corollary}
\textbf{Proof.} Substitute  $(q-1)x$ instead of $x$ in (\ref{eq:2phi1-identity}) and divide both sides of by $(q-1)^{\alpha}$.
Now taking the limit as $q\to1^{-}$ yields (\ref{eq:1F1-identity}).$\hfill\square$

\section{Proofs}
We will repeatedly use the identity
\begin{equation}\label{eq:q-poch-gamma}
\frac{\Gamma_q(x+k)}{\Gamma_q(x)}=\frac{(q^x;q)_k}{(1-q)^k}
\end{equation}
valid for any real $x\ne0,-1,-2,\ldots$ and $k\in\N_0$ which is easy to verify with the help of the definition (\ref{eq:q-Gamma-def}).

In order to give a proof of Theorem~\ref{th:tphis} we need to do some preliminary work.
Define the rational function
\begin{equation}\label{eq:ratf}
R_{t,s}(y)=\frac{\prod_{k=1}^t(c_k+y)}{\prod_{k=1}^s(d_k+y)}
\end{equation}
with nonnegative parameters $c_k$ and $d_k$.  Recall that elementary symmetric polynomials $e_k(\c)=e_k(c_1,\dots,c_r)$, $k=0,1,\ldots,r$, have been defined in (\ref{eq:symmetricpol}).  We will need the following lemma.
\begin{lemma}\label{l:inrdecr}\emph{(\cite[Lemma~3]{KKJA2016})}
If $t\geq s$ and  conditions \emph{(\ref{eq:ratincr})} are satisfied, then the function $R_{t,s}(y)$ is monotone increasing on $(0,\infty)$.

If $t\leq s$ and  conditions \emph{(\ref{eq:ratdecr})} are satisfied, then the function $R_{t,s}(y)$ is monotone decreasing on $(0,\infty)$.
\end{lemma}

The two propositions below can be found in \cite[Lemmas 2 and 3]{KK2}.
\begin{lemma}\label{l:disr1}
Let $f$ be a nonnegative-valued function defined on $[0,\infty)$ and
$$
\Delta_f(\alpha,\beta)=f(\mu+\alpha)f(\mu+\beta)-f(\mu)f(\mu+\beta+\alpha)\geq 0  \ \ \text{for} \ \ \alpha=1 \ \ \text{and all} \ \ \mu,\beta\geq 0.
$$
Then $\Delta_f(\alpha,\beta)\geq0$ for all $\alpha\in\N$ and $\mu,\beta\geq 0$. If inequality is strict in the hypothesis of the lemma then it is also strict in the conclusion.
\end{lemma}

\begin{lemma}\label{l:disr2} Let $f$ be defined by the formal series
$$
f(\mu;x)=\sum_{k=0}^{\infty} f_k(\mu)x^k,
$$
where the coefficients $f_k(\mu)$ are continuous nonnegative functions.  Assume that $\Delta_f(1,\beta;x)$ defined in \emph{(\ref{eq:delta-defined})} has nonnegative \emph{(}non-positive\emph{)} coefficients at all powers of $x$ for $\mu,\beta\ge0$. Then $\Delta_f(\alpha,\beta;x)$ has nonnegative \emph{(}non-positive\emph{)} coefficients at powers of $x$ for all $\alpha\in\N$,
$\alpha\le\beta+1$ and $\mu\ge0$.
\end{lemma}

\noindent\textbf{Proof of Theorem~\ref{th:tphis}.}
According to Lemmas~\ref{l:disr1} and \ref{l:disr2} it suffices to consider the case $\alpha=1$, $\beta\ge0$.
We have:
\begin{multline*}
\frac{\Gamma_q(\b+\mu)\Gamma_q(\b+\mu+\beta)}{\Gamma_q(\a+\mu)\Gamma_q(\a+\mu+\beta)(1-q)^{s-t}}\Delta_{g}(1,\beta;x)
\\
= \frac{(1-q^{\a+\mu})}{(1-q^{\b+\mu})}
{}_t\phi_{s}\left(\!\begin{array}{c}q^{\a+\mu+1}\\q^{\b+\mu+1}\!\end{array}\vline\,(q-1)^{1+s-t}x\right) {}_t\phi_{s}\left(\!\begin{array}{c}q^{\a+\mu+\beta}\\q^{\b+\mu+\beta}\!\end{array}\vline\,(q-1)^{1+s-t}x\right)
\\
-\frac{(1-q^{\a+\mu+\beta})}{(1-q^{\b+\mu+\beta})}
{}_t\phi_{s}\left(\!\begin{array}{c}q^{\a+\mu}\\q^{\b+\mu}\!\end{array}\vline\,(q-1)^{1+s-t}x\right) {}_t\phi_{s}\left(\!\begin{array}{c}q^{\a+\mu+\beta+1}\\q^{\b+\mu+\beta+1}\!\end{array}\vline\,(q-1)^{1+s-t}x\right)
\\
=\sum_{j=1}^{\infty}\frac{(q^{\a+\mu};q)_j}{(q^{\b+\mu};q)_j}\left[q^{\binom{j-1}{2}}\right]^{1+s-t}
\frac{((1-q)^{1+s-t}x)^{j-1}}{(q;q)_{j-1}} \sum_{j=0}^{\infty}\frac{(q^{\a+\mu+\beta};q)_j}{(q^{\b+\mu+\beta};q)_j}\left[q^{\binom{j}{2}}\right]^{1+s-t}\frac{((1-q)^{1+s-t}x)^{j}}{(q;q)_{j}} \\
-
\sum_{j=0}^{\infty}\frac{(q^{\a+\mu};q)_j}{(q^{\b+\mu};q)_j}\left[q^{\binom{j}{2}}\right]^{1+s-t}\frac{((1-q)^{1+s-t}x)^{j}}{(q;q)_{j}}
\sum_{j=1}^{\infty}\frac{(q^{\a+\mu+\beta};q)_j}{(q^{\b+\mu+\beta};q)_j}\left[q^{\binom{j-1}{2}}\right]^{1+s-t}\frac{((1-q)^{1+s-t}x)^{j-1}}{(q;q)_{j-1}}
\\
=\sum_{m=1}^{\infty}((1-q)^{1+s-t}x)^{m-1}\sum_{k=0}^{m} \frac{(q^{\a+\mu};q)_{k}(q^{\a+\mu+\beta};q)_{m-k}(1-q^k)}{(q^{\b+\mu};q)_{k}(q^{\b+\mu+\beta};q)_{m-k}(q;q)_k(q;q)_{m-k}}\left[q^{\left(\binom{k-1}{2}+\binom{m-k}{2}\right)}\right]^{1+s-t}
\\
-\sum_{m=1}^{\infty}
((1-q)^{1+s-t}x)^{m-1}\sum_{k=0}^{m}\frac{(q^{\a+\mu};q)_k(q^{\a+\mu+\beta};q)_{m-k}(1-q^{m-k})}{(q^{\b+\mu};q)_k(q^{\b+\mu+\beta};q)_{m-k}(q;q)_k(q;q)_{m-k}}\left[q^{\left(\binom{k}{2}+\binom{m-k-1}{2}\right)}\right]^{1+s-t}
\\
=\sum_{m=1}^{\infty}((1-q)^{1+s-t}x)^{m-1}\sum_{k=0}^{m}\frac{(q^{\a+\mu};q)_k(q^{\a+\mu+\beta};q)_{m-k}}{(q^{\b+\mu};q)_k(q^{\b+\mu+\beta};q)_{m-k}(q;q)_k(q;q)_{m-k}}
\\
\times\left\{(1-q^k)\left[q^{\left(\binom{k-1}{2}+\binom{m-k}{2}\right)}\right]^{1+s-t}-(1-q^{m-k})\left[q^{\left(\binom{k}{2}+\binom{m-k-1}{2}\right)}\right]^{1+s-t}\right\}.
\end{multline*}
Applying the Gauss pairing to the finite inner sum, the last expression becomes
\begin{multline*}
\sum_{m=1}^{\infty}((1-q)^{1+s-t}x)^{m-1}
\sum_{0\leq k\leq m/2}\frac{1}{{(q;q)_k(q;q)_{m-k}}}\times
\\
\left\{(1-q^{m-k})\left[q^{\left(\binom{k}{2}+\binom{m-k-1}{2}\right)}\right]^{1+s-t}-(1-q^k)\left[q^{\left(\binom{k-1}{2}+\binom{m-k}{2}\right)}\right]^{1+s-t}\right\}
\\
\times \left[\frac{(q^{\a+\mu};q)_{m-k}(q^{\a+\mu+\beta};q)_{k}}{(q^{\b+\mu};q)_{m-k}(q^{\b+\mu+\beta};q)_{k}}-\frac{(q^{\a+\mu};q)_k(q^{\a+\mu+\beta};q)_{m-k}}{(q^{\b+\mu};q)_k(q^{\b+\mu+\beta};q)_{m-k}}\right].
\end{multline*}
The apparent (for even $m$) unpaired middle term  vanishes due to the factor
$$
(1-q^{m-k})\left[q^{\left(\binom{k}{2}+\binom{m-k-1}{2}\right)}\right]^{1+s-t}-(1-q^k)\left[q^{\left(\binom{k-1}{2}+\binom{m-k}{2}\right)}\right]^{1+s-t},
$$
which is equal to zero for $k=m-k$ and is non-negative for $k\leq m-k$, $1+s-t\ge0$ because, clearly, $1-q^{m-k}\ge1-q^k$ and
\begin{multline*}
\left[q^{\left(\binom{k}{2}+\binom{m-k-1}{2}\right)}\right]^{1+s-t} \ge \left[q^{\left(\binom{k-1}{2}+\binom{m-k}{2}\right)}\right]^{1+s-t}
~\Leftrightarrow~q^{\left(\binom{k}{2}+\binom{m-k-1}{2}\right)} \ge q^{\left(\binom{k-1}{2}+\binom{m-k}{2}\right)}
\\
\Leftrightarrow~\binom{k}{2}+\binom{m-k-1}{2} \le \binom{k-1}{2}+\binom{m-k}{2}~\Leftrightarrow~2(2k-m)\le0.
\end{multline*}
It remains to show that the factor in brackets has a constant sign. Indeed, for $k\le{m-k}$ we can rearrange this factor as follows:
\begin{multline*}
\frac{(q^{\a+\mu};q)_{m-k}(q^{\a+\mu+\beta};q)_{k}}{(q^{\b+\mu};q)_{m-k}(q^{\b+\mu+\beta};q)_{k}}
-\frac{(q^{\a+\mu};q)_k(q^{\a+\mu+\beta};q)_{m-k}}{(q^{\b+\mu};q)_k(q^{\b+\mu+\beta};q)_{m-k}}
\\
=\frac{(q^{\a+\mu};q)_k(q^{\a+\mu+\beta};q)_k}{(q^{\b+\mu};q)_k(q^{\b+\mu+\beta};q)_k}
\biggl\{\prod\nolimits_{j=k}^{m-k-1}h(q^{\mu+j})-\prod\nolimits_{j=k}^{m-k-1}h(q^{\mu+\beta+j})\biggr\}.
\end{multline*}
Here the function $h(z)$ is defined for $0<z<1$  by
\begin{multline*}
h(z)=\frac{(1-zq^{\a})}{(1-zq^{\b})}=\frac{\prod_{j=1}^{t}(1-z q^{a_j})}{\prod_{j=1}^{s}(1-z q^{b_j})}=
\frac{\prod_{j=1}^{t}q^{a_j}\prod_{j=1}^{t}(q^{-a_j}-z)}{\prod_{j=1}^{s}q^{b_j}\prod_{j=1}^{s}(q^{-b_j}-z)}
\\
=\frac{\prod_{j=1}^{t}q^{a_j}\prod_{j=1}^{t}((q^{-a_j}-1)+(1-z))}{\prod_{j=1}^{s}q^{b_j}\prod_{j=1}^{s}((q^{-b_j}-1)+(1-z))}
=\frac{\prod_{j=1}^{t}q^{a_j}\prod_{j=1}^{t}(c_j+y)}{\prod_{j=1}^{s}q^{b_j}\prod_{j=1}^{s}(d_j+y)}=\frac{q^{\a}}{q^{\b}}R_{t,s}(y),
\end{multline*}
where $c_j=q^{-a_j}-1>0$, $d_j=q^{-b_j}-1>0$, $y=1-z>0$ and $R_{t,s}$ is given in (\ref{eq:ratf}). Note, that the function $y(z)$ is decreasing.  The claims (a) and (b) now follow by Lemma~\ref{l:inrdecr} in view of $\mu,\beta\ge0$.
$\hfill\square$

\bigskip

\noindent\textbf{Proof of Theorem~\ref{th:1phi1}.} Indeed, if
$$
\Delta_f(\alpha,\beta;x)=f(\mu+\alpha;x)f(\mu+\beta;x)-f(\mu;x)f(\mu+\alpha+\beta;x)=\sum\limits_{m=0}^{\infty}\delta_mx^m,
$$
then, by  using the Cauchy product and Lemma~\ref{cr:qBessel},
\begin{multline*}
\delta_m=\sum\limits_{k=0}^{m}\frac{1}{(q;q)_k(q;q)_{m-k}}
\left\{\frac{1}{(q^{\mu+\alpha};q)_k(q^{\mu+\beta};q)_{m-k}}-\frac{1}{(q^{\mu};q)_k(q^{\mu+\alpha+\beta};q)_{m-k}}\right\}
\\
=\sum\limits_{k=0}^{m}\frac{(q^{(2\mu+\alpha+\beta-1)/2},q^{(2\mu+\alpha+\beta)/2},-q^{(2\mu+\alpha+\beta-1)/2},-q^{(2\mu+\alpha+\beta)/2};q)_k}
{(q^{2\mu+\alpha+\beta-1};q)_k(q;q)_{k}(q;q)_{m-k}}\times
\\
\left\{\frac{1}{(q^{\mu+\alpha},q^{\mu+\beta};q)_k}-\frac{1}{(q^{\mu},q^{\mu+\alpha+\beta};q)_k}\right\}.
\end{multline*}
It remains to note that according to definition (\ref{eq:q-shifted}) the inequality
$$
\frac{1}{(q^{\mu+\alpha},q^{\mu+\beta};q)_k}<\frac{1}{(q^{\mu},q^{\mu+\alpha+\beta};q)_k}
$$
for $\mu,\alpha,\beta>0$ is implied by
$$
(1-q^{\mu+j})(1-q^{\mu+\alpha+\beta+j})<(1-q^{\mu+\alpha+j})(1-q^{\mu+\beta+j})~\Leftrightarrow~q^{\alpha}+q^{\beta}<1+q^{\alpha+\beta}
$$
for $j=0,1,\ldots,k-1$.  The last inequality is straightforward.
$\hfill\square$

\bigskip

\noindent\textbf{Proof of Theorem~\ref{th:2phi1logconcavity}.} We have
$$
\Delta_{\tilde{f}}(\alpha,\beta;x)=\tilde{f}(\mu+\alpha;x)\tilde{f}(\mu+\beta;x)-\tilde{f}(\mu;x)\tilde{f}(\mu+\alpha+\beta;x)
=\sum\limits_{m=0}^{\infty}[x/(1-q)]^m\tilde{\delta}_m.
$$
Using the Cauchy product and Lemma~\ref{cr:qBessel} the coefficients $\tilde{\delta}_m$ are computed as follows:
\begin{multline*}
\tilde{\delta}_m=
\sum\limits_{k=0}^{m}\frac{1}{(q;q)_k(q;q)_{m-k}}\times
\\
\left\{\frac{1}{\Gamma_q(\mu+\alpha+k)\Gamma_q(\mu+\beta+m-k)}
-\frac{1}{\Gamma_q(\mu+k)\Gamma_q(\mu+\alpha+\beta+m-k)}
\right\}
\\=
\sum\limits_{k=0}^{m}\frac{(1-q)^m}{(q;q)_k(q;q)_{m-k}}
\left\{\frac{[\Gamma_q(\mu+\alpha)\Gamma_q(\mu+\beta)]^{-1}}{(q^{\mu+\alpha};q)_k(q^{\mu+\beta};q)_{m-k}}
-\frac{[\Gamma_q(\mu)\Gamma_q(\mu+\alpha+\beta)]^{-1}}{(q^{\mu};q)_k(q^{\mu+\alpha+\beta};q)_{m-k}}\right\}
\\
=\sum\limits_{k=0}^{m}\frac{(1-q)^m(q^{(2\mu+\alpha+\beta-1)/2},q^{(2\mu+\alpha+\beta)/2},-q^{(2\mu+\alpha+\beta-1)/2},-q^{(2\mu+\alpha+\beta)/2};q)_k}
{(q^{2\mu+\alpha+\beta-1};q)_k(q;q)_{k}(q;q)_{m-k}}\times
\\
\biggl\{\frac{[\Gamma_q(\mu+\alpha)\Gamma_q(\mu+\beta)]^{-1}}{(q^{\mu+\alpha},q^{\mu+\beta};q)_k}
-\frac{[\Gamma_q(\mu)\Gamma_q(\mu+\alpha+\beta)]^{-1}}{(q^{\mu},q^{\mu+\alpha+\beta};q)_k}\biggr\}
\\
=\frac{(1-q)^{2\mu+\alpha+\beta-2}}{[(q;q)_{\infty}]^2}\sum\limits_{k=0}^{m}
\frac{(1-q)^m(q^{(2\mu+\alpha+\beta-1)/2},q^{(2\mu+\alpha+\beta)/2},-q^{(2\mu+\alpha+\beta-1)/2},-q^{(2\mu+\alpha+\beta)/2};q)_k}
{(q^{2\mu+\alpha+\beta-1};q)_k(q;q)_{k}(q;q)_{m-k}}\times
\\
\biggl\{\frac{(q^{\mu+\alpha},q^{\mu+\beta};q)_{\infty}}{(q^{\mu+\alpha},q^{\mu+\beta};q)_k}
-\frac{(q^{\mu},q^{\mu+\alpha+\beta};q)_{\infty}}{(q^{\mu},q^{\mu+\alpha+\beta};q)_k}\biggr\},
\end{multline*}
where we applied (\ref{eq:q-poch-gamma}) in the first equality and  (\ref{eq:qBesselfinitesum}) in the second.
The expression in braces on the right hand side is immediately seen to reduce to
$$
(q^{\mu+\alpha+k},q^{\mu+\beta+k};q)_{\infty}-(q^{\mu+k},q^{\mu+\alpha+\beta+k};q)_{\infty}.
$$
To prove its positivity, it suffices to show that each factor in the first term is greater than that in the second.
Indeed, for any $j\in\N_0$
$$
(1-q^{\mu+\alpha+k+j})(1-q^{\mu+\beta+k+j})>(1-q^{\mu+k+j})(1-q^{\mu+\alpha+\beta+k+j}),
$$
which is seen by expanding both sides and applying the elementary inequality $u+v<1+uv$ valid for $u,v\in(0,1)$.
$\hfill\square$

\bigskip

The summation formula contained in the following lemma was established in \cite[Lemma~8]{KKJMAA2018}.
\begin{lemma}\label{l:recqgamma}
The following identity holds true\emph{:}
\begin{multline}\label{eq:recqgamma}
\sum_{k=0}^m\left(\frac{1}{\Gamma_q(k+\mu+1)\Gamma_q(m-k+\mu+\beta)}-\frac{1}{\Gamma_q(k+\mu)\Gamma_q(m-k+\mu+\beta+1)}\right)
\\
=\frac{\left(q^{\mu+\beta};q\right)_{m+1}-\left(q^{\mu};q\right)_{m+1}}{\Gamma_q(\mu+m+1)\Gamma_q(\mu+\beta+m+1)(1-q)^{m+1}}.
\end{multline}
\end{lemma}

\noindent\textbf{Proof of Theorem~\ref{th:2phi1-identity}.}
First, consider the case $\alpha=1$. Using (\ref{eq:q-poch-gamma}) we  have:
\begin{multline*}
\frac{1}{\Gamma_q(\mu+1)\Gamma_q(\mu+\beta)}
{}_2\phi_1\left(\!\begin{array}{c}q, 0\\q^{\mu+1}\!\end{array}\vline\,x\right)
{}_2\phi_1\left(\!\begin{array}{c}q, 0\\q^{\mu+\beta}\!\end{array}\vline\,x\right)
\\
-\frac{1}{\Gamma_q(\mu)\Gamma_q(\mu+1+\beta)}
{}_2\phi_1\left(\!\begin{array}{c}q, 0\\q^{\mu}\!\end{array}\vline\,x\right)
{}_2\phi_1\left(\!\begin{array}{c}q, 0\\q^{\mu+1+\beta}\!\end{array}\vline\,x\right)=
\\
\frac{1}{\Gamma_q(\mu+1)\Gamma_q(\mu+\beta)} \sum_{k=0}^{\infty}\frac{x^k}{\left(q^{\mu+1};q\right)_k}\sum_{n=0}^{\infty}\frac{x^n}{\left(q^{\mu+\beta};q\right)_n}
-\frac{1}{\Gamma_q(\mu)\Gamma_q(\mu+1+\beta)}\sum_{k=0}^{\infty}\frac{x^k}{\left(q^{\mu};q\right)_k}\sum_{n=0}^{\infty}\frac{x^n}{(q^{\mu+1+\beta};q)_n}
\\
=\sum_{m=0}^{\infty}\frac{x^m}{(1-q)^m}\sum_{k=0}^m\left\{\frac{(1-q)^m}{\left(q^{1+\mu};q\right)_{k}\Gamma_q(1+\mu)\left(q^{\mu+\beta};q\right)_{m-k}\Gamma_q(\mu+\beta)} \right.
\\-\left.\frac{(1-q)^m}{\left(q^{\mu};q\right)_k\Gamma_q(\mu)\left(q^{\mu+1+\beta};q\right)_{m-k}\Gamma_q(\mu+1+\beta)}\right\}
\\
=\sum_{m=0}^{\infty}\frac{x^m}{(1-q)^m}\sum_{k=0}^m\left\{\frac{1}{\Gamma_q(k+\mu+1)\Gamma_q(m-k+\mu+\beta)}
-\frac{1}{\Gamma_q(k+\mu)\Gamma_q(m-k+\mu+\beta+1)}\right\}.
\end{multline*}
We can now apply Lemma~\ref{l:recqgamma} to the expression on the right hand side to get
\begin{multline*}
\sum_{m=0}^{\infty}\frac{x^m}{(1-q)^{2m+1}}\left\{\frac{\left(q^{\mu+\beta};q\right)_{m+1}-\left(q^{\mu};q\right)_{m+1}}{\Gamma_q(\mu+m+1)\Gamma_q(\mu+\beta+m+1)}\right\}
\\
=\frac{1}{\Gamma_q(\mu+1)\Gamma_q(\mu+\beta)}\sum_{m=0}^{\infty}\frac{\left(q;q\right)_{m}}{\left(q^{\mu+1},q;q\right)_m}x^m
-\frac{1}{\Gamma_q(\mu)\Gamma_q(\mu+1+\beta)}\sum_{m=0}^{\infty}\frac{\left(q;q\right)_{m}}{\left(q^{\mu+\beta+1},q;q\right)_m} x^m
\\
=\frac{1}{\Gamma_q(\mu+1)\Gamma_q(\mu+\beta)}
{}_2\phi_1\left(\!\begin{array}{c}q, 0\\q^{\mu+1}\!\end{array}\vline\,x\right)-\frac{1}{\Gamma_q(\mu)\Gamma_q(\mu+1+\beta)}
{}_2\phi_1\left(\!\begin{array}{c}q, 0\\q^{\mu+\beta+1}\!\end{array}\vline\,x\right),
\end{multline*}
where (\ref{eq:q-poch-gamma}) has been applied in the first equality. Hence, the theorem is proved for $\alpha=1$.

Next, define the function
$$
h(\mu):=\frac{1}{\Gamma_q(\mu)}{}_2\phi_1\left(\!\begin{array}{c}q, 0\\q^{\mu}\!\end{array}\vline\,x\right).
$$
In terms of this function we compute the generalized Tur\'{a}nian as follows:
\begin{multline*}
h(\mu+\alpha)h(\mu+\beta)-h(\mu)h(\mu+\alpha+\beta)
\\
=[h(\mu+\alpha)h(\mu+\beta)-h(\mu+\alpha-1)h(\mu+\beta+1)]
\\
+[h(\mu+\alpha-1)h(\mu+\beta+1)-h(\mu+\alpha-2)h(\mu+\beta+2)]+
\\
\cdots+[h(\mu+2)h(\mu+\beta+\alpha-2)-h(\mu+1)h(\mu+\alpha+\beta-1)]
\\
+[h(\mu+1)h(\mu+\beta+\alpha-1)-h(\mu)h(\mu+\alpha+\beta)]
\\
=\sum_{j=0}^{\alpha-1}\biggl\{\frac{1}{\Gamma_q(\mu+1+j)\Gamma_q(\mu+\alpha+\beta-1-j)}
{}_2\phi_1\left(\!\begin{array}{c}q, 0\\q^{\mu+1+j}\!\end{array}\vline\,x\right)
\\
-\frac{1}{\Gamma_q(\mu+j)\Gamma_q(\mu+\alpha+\beta-j)}
{}_2\phi_1\left(\!\begin{array}{c}q, 0\\q^{\mu+\alpha+\beta-j}\!\end{array}\vline\,x\right)\biggr\},
\end{multline*}
where we applied the $\alpha=1$ case established above to each bracketed expression. It remains to multiply throughout  by $\Gamma_q(\mu+\alpha)\Gamma_q(\mu+\alpha+\beta)$ to complete the proof of the theorem.
$\hfill\square$

\bigskip

{\bf Acknowledgements.} This research has been supported by the Russian Science Foundation under project 14-11-00022.

\end{document}